\documentclass[english,11pt]{article}
\usepackage[latin1]{inputenc}
\usepackage{amsmath,amsthm,amssymb,babel}

\newtheorem{teo}{Theorem}

\newtheorem{lemma}{Lemma}

\def\proof{{\it Proof.}\ }
\def\endproof{\hfill $\Box$\par\vskip3mm}

\def\eq#1{(\ref{#1})}

\def\neweq#1{\begin{equation}\label{#1}}
\def\endeq{\end{equation}}

\def\phi{\varphi}
\def\RR{{\mathbb R} }

\def\di{\displaystyle}

\date{}

\title{\sc On a nonhomogeneous quasilinear eigenvalue problem in Sobolev spaces
with variable exponent\thanks{
Correspondence address: Vicen\c{t}iu R\u{a}dulescu, Department of
Mathematics, University of Craiova,  200585 Craiova, Romania. E-mail: 
{\tt
vicentiu.radulescu@math.cnrs.fr}}}
\author{\sc Mihai Mih\u ailescu and Vicen\c{t}iu R\u{a}dulescu\\
\small
Department of Mathematics, University of Craiova,  200585 Craiova, 
Romania\\
\small
E-mail addresses: {\tt mmihailes@yahoo.com}\qquad {\tt 
vicentiu.radulescu@math.cnrs.fr}}

\begin{document}
\baselineskip16pt
\maketitle
\noindent{\small{\sc Abstract}. We consider the nonlinear eigenvalue problem 
$-{\rm div}\left(|\nabla u|^{p(x)-2}\nabla u\right)=\lambda |u|^{q(x)-2}u$ in $\Omega$,
$u=0$ on $\partial\Omega$, where $\Omega$ is a bounded open set in $\RR^N$
with smooth boundary and $p$, $q$ are continuous functions on $\overline\Omega$ such that $1<\inf_\Omega q<
\inf_\Omega p<\sup_\Omega q$, $\sup_\Omega p<N$, and $q(x)<Np(x)/\left(N-p(x)\right)$ for all $x\in\overline\Omega$. 
The main result of this paper establishes that any $\lambda>0$ sufficiently small is an eigenvalue of the above nonhomogeneous
quasilinear problem. The proof relies on simple variational arguments based on Ekeland's variational principle.
}

\section{Introduction and preliminary results}
 
A basic result in the elementary theory of linear partial differential equations asserts that the spectrum of the Laplace
operator in $H^1_0(\Omega)$ is discrete, where $\Omega$ is a bounded open set in $\RR^N$ with smooth boundary. More precisely, the
problem
$$
\left\{\begin{array}{lll}
-\Delta u=\lambda u\qquad &\mbox{in}& \Omega\\
u=0\qquad &\mbox{on}& \partial\Omega
\end{array}\right.
$$
has an unbounded sequence of eigenvalues $0<\lambda_1<\lambda_2\leq\ldots\leq\lambda_n\leq\ldots$. This celebrated result
goes back to the Riesz-Fredholm theory of self-adjoint and compact operators on Hilbert spaces. The anisotropic case
$$
\left\{\begin{array}{lll}
-\Delta u=\lambda a(x) u\qquad &\mbox{in}& \Omega\\
u=0\qquad &\mbox{on}& \partial\Omega
\end{array}\right.
$$
was considered by Bocher \cite{boc}, Hess and Kato \cite{hes}, Minakshisundaram and Pleijel \cite{minple,ple}.
For instance, Minakshisundaram and  Pleijel proved that the above eigenvalue problem has an unbounded sequence of 
positive eigenvalues if $a\in L^\infty (\Omega)$, $a\geq 0$ in $\Omega$, 
and $a>0$ in $\Omega_0\subset\Omega$, where $|\Omega_0|>0$. Eigenvalue problems for homogeneous quasilinear problems have been
intensively studied in the last decades (see, e.g., Anane \cite{anane}). 
 
This paper is motivated by recent advances in elastic mechanics and electrorheological fluids 
(sometimes referred to as ``smart fluids"), where some processes are modeled by nonhomogeneous quasilinear operators (see
Diening \cite{D},
Halsey \cite{hal}, Ruzicka \cite{R}, Zhikov \cite{Z1}, and the references therein). We refer mainly to the
 $p(x)$--Laplace operator 
$\Delta_{p(x)}u:={\rm div}(|\nabla u|^{p(x)-2}\nabla u)$, where $p$ is a continuous non-constant function.
This differential operator is a natural generalization of the 
$p$-Laplace operator $\Delta_p u:={\rm div}\,(|\nabla u|^{p-2}\nabla u)$, where
$p>1$ is a real constant. However, the $p(x)$-Laplace operator possesses more complicated nonlinearities 
than the $p$-Laplace operator, due to the fact that $\Delta_{p(x)}$ is not homogeneous. 
 
In this paper we are concerned with the nonhomogeneous eigenvalue problem
\begin{equation}\label{2}
\left\{\begin{array}{lll}
-{\rm div}(|\nabla u|^{p(x)-2}\nabla u)=\lambda|u|^{q(x)-2}u, &\mbox{for}&
x\in\Omega\\
u=0, &\mbox{for}& x\in\partial\Omega\,,
\end{array}\right.
\end{equation}
where $\Omega\subset\RR^N$ ($N\geq 3$) is a bounded domain with
smooth boundary, $\lambda>0$ is a real number, and $p$, $q$ are continuous on $\overline\Omega$. 

The case $p(x)=q(x)$ was considered by Fan, Zhang and Zhao in \cite{FZZ} who, using
the Ljusternik-Schnirelmann critical point theory, established the existence of 
a sequence of eigenvalues. 
Denoting by $\Lambda$ the set of all 
nonnegative eigenvalues, Fan, Zhang and Zhao showed that $\sup\Lambda=+\infty$ and they pointed out that 
only under additional assumptions we have
$\inf\Lambda>0$. We remark that for the $p$-Laplace 
operator (corresponding to $p(x)\equiv p$) we always have $\inf\Lambda>0$.

In this paper we  study problem \eq{2} under the basic assumption
\begin{equation}\label{3}
1<\min_{x\in\overline\Omega}q(x)<\min_{x\in\overline\Omega}p(x)<\max_{x\in\overline\Omega}q(x).
\end{equation}
Our main result establishes the existence of a continuous family of eigenvalues for problem \eq{2} 
in a neighborhood of the origin. More precisely, we show that there exists $\lambda^\star>0$ 
such that {\it any} $\lambda\in(0,\lambda^\star)$ is an eigenvalue for problem \eq{2}.
    
We start with some preliminary basic results on the theory of
 Lebesgue--Sobolev spaces with variable exponent. For more details we refer to the book by Musielak 
\cite{M} and the papers by Edmunds et al. \cite{edm, edm2, edm3}, Kovacik and  R\'akosn\'{\i}k \cite{KR},
Mih\u ailescu and R\u adulescu \cite{RoyalSoc}, and Samko and Vakulov \cite{samko}. 

Assume that $p\in C(\overline\Omega)$ and $p(x)>1$, for all $x\in\overline\Omega$.

Set
$$C_+(\overline\Omega)=\{h;\;h\in C(\overline\Omega),\;h(x)>1\;{\rm 
for}\;
{\rm all}\;x\in\overline\Omega\}.$$
For any $h\in C_+(\overline\Omega)$ we define
$$h^+=\sup_{x\in\Omega}h(x)\qquad\mbox{and}\qquad h^-=
\inf_{x\in\Omega}h(x).$$
For any $p(x)\in C_+(\overline\Omega)$, we define the variable exponent
Lebesgue space
$$L^{p(x)}(\Omega)=\{u;\ u\ \mbox{is a
 measurable real-valued function such that }
\int_\Omega|u(x)|^{p(x)}\;dx<\infty\}.$$
We define a norm, the so-called {\it Luxemburg norm}, on this space by 
the
formula
$$|u|_{p(x)}=\inf\left\{\mu>0;\;\int_\Omega\left|
\frac{u(x)}{\mu}\right|^{p(x)}\;dx\leq 1\right\}.$$
We remember that the variable exponent Lebesgue spaces are separable and 
reflexive Banach spaces. If $0 <|\Omega|<\infty$ and $p_1$, $p_2$
are variable exponent so that $p_1(x) \leq p_2(x)$ almost everywhere in
$\Omega$ then there exists the continuous embedding
$L^{p_2(x)}(\Omega)\hookrightarrow L^{p_1(x)}(\Omega)$.

We denote by $L^{p^{'}(x)}(\Omega)$ the conjugate space
of $L^{p(x)}(\Omega)$, where $1/p(x)+1/p^{'}(x)=1$. For any
$u\in L^{p(x)}(\Omega)$ and $v\in L^{p^{'}(x)}(\Omega)$ the H\"older
type inequality
\begin{equation}\label{Hol}
\left|\int_\Omega uv\;dx\right|\leq\left(\frac{1}{p^-}+
\frac{1}{{p^{'}}^-}\right)|u|_{p(x)}|v|_{p^{'}(x)}
\end{equation}
holds true.

An important role in manipulating the generalized Lebesgue-Sobolev 
spaces
is played by the {\it modular} of the $L^{p(x)}(\Omega)$ space, which 
is
the mapping
 $\rho_{p(x)}:L^{p(x)}(\Omega)\rightarrow\RR$ defined by
$$\rho_{p(x)}(u)=\int_\Omega|u|^{p(x)}\;dx.$$
If $(u_n)$, $u\in L^{p(x)}(\Omega)$ then the following relations
hold true
\begin{equation}\label{L4}
|u|_{p(x)}>1\;\;\;\Rightarrow\;\;\;|u|_{p(x)}^{p^-}\leq\rho_{p(x)}(u)
\leq|u|_{p(x)}^{p^+}
\end{equation}
\begin{equation}\label{L5}
|u|_{p(x)}<1\;\;\;\Rightarrow\;\;\;|u|_{p(x)}^{p^+}\leq
\rho_{p(x)}(u)\leq|u|_{p(x)}^{p^-}
\end{equation}
\begin{equation}\label{L6}
|u_n-u|_{p(x)}\rightarrow 0\;\;\;\Leftrightarrow\;\;\;\rho_{p(x)}
(u_n-u)\rightarrow 0.
\end{equation}

Next, we define $W_0^{1,p(x)}(\Omega)$ as the closure of
$C_0^\infty(\Omega)$ under the norm
$$\|u\|=|\nabla u|_{p(x)}.$$
The space $(W_0^{1,p(x)}(\Omega),\|\cdot\|)$ is a separable and
reflexive Banach space. We note that if $s(x)\in C_+(\overline\Omega)$
and $s(x)<p^\star(x)$ for all $x\in\overline\Omega$ then the
embedding
$W_0^{1,p(x)}(\Omega)\hookrightarrow L^{s(x)}(\Omega)$
is compact and continuous, where $p^\star(x)=\frac{Np(x)}{N-p(x)}$
if $p(x)<N$ or $p^\star(x)=+\infty$ if $p(x)\geq N$.

\section{The main result}
We say that $\lambda\in\RR$ is an eigenvalue of problem \eq{2} if 
there exists $u\in W_0^{1,p(x)}(\Omega)
\setminus\{0\}$ such that
$$\int_{\Omega}|\nabla u|^{p(x)-2}\nabla u\nabla v\;dx-\lambda\int_{\Omega} 
|u|^{q(x)-2}uv\;dx=0,$$
for all $v\in W_0^{1,p(x)}(\Omega)$. We point out that if $\lambda$ is an 
eigenvalue of the problem \eq{2} then the corresponding $u\in W_0^{1,p(x)}(\Omega)\setminus\{0\}$ 
is a  weak solution of \eq{2}. 

Our main result is given by the following theorem.

\begin{teo}\label{t1}
Assume that condition \eq{3} is fulfilled, $\max_{x\in\overline\Omega}p(x)<N$ and 
$q(x)<p^\star(x)$ for all $x\in\overline\Omega$. Then there exists $\lambda^\star>0$ 
such that any $\lambda\in(0,\lambda^\star)$ is an eigenvalue for problem \eq{2}.
\end{teo}

The above result implies
$$\inf_{u\in W_0^{1,p(x)}(\Omega)\setminus\{0\}}\frac{\di\int_\Omega
|\nabla u|^{p(x)}\;dx}{\di\int_\Omega|u|^{q(x)}\;dx}=0.$$
Thus,  for any positive constant $C$, there exists $u_0\in 
W_0^{1,p(x)}(\Omega)$ such that 
$$C\int_\Omega|u_0|^{q(x)}\;dx\geq\int_\Omega|\nabla u_0|^{p(x)}\;dx.$$

\medskip
Let $E$ denote the generalized Sobolev space $W_0^{1,p(x)}(\Omega)$.

For any $\lambda>0$ the energy functional corresponding to problem \eq{2} is defined as
$J_\lambda:E\rightarrow\RR$,
$$J_\lambda(u)=\int_{\Omega}\frac{1}{p(x)}|\nabla u|^{p(x)}\;dx-\lambda
\int_\Omega\frac{1}{q(x)}|u|^{q(x)}\;dx.$$
Standard arguments imply that 
$J_\lambda\in C^1(E,\RR)$ and
$$\langle J_\lambda^{'}(u),v\rangle=\int_{\Omega}|\nabla u|^{p(x)-2}
\nabla u\nabla v\;dx-\lambda\int_\Omega|u|^{q(x)-2}uv\;dx,$$
for all $u,\;v\in E$. Thus the weak solutions of \eq{2} coincide with the critical 
points of $J_\lambda$. If such an weak solution exists and is nontrivial then 
the corresponding $\lambda$ is an eigenvalue of problem \eq{2}.

\begin{lemma}\label{l1}
There exists $\lambda^\star>0$ such that for any $\lambda\in(0,\lambda^\star)$ 
there exist $\rho$, $a>0$  such that $J_\lambda(u)\geq a>0$ 
for any $u\in E$ with $\|u\|=\rho$.
\end{lemma}
\proof
Since $q(x)<p^\star(x)$ for all $x\in\overline\Omega$ it follows that $E$ is continuously 
embedded in $L^{q(x)}(\Omega)$. So, there exists a positive constant $c_1$ such that
\begin{equation}\label{el1}
|u|_{q(x)}\leq c_1\|u\|,\;\;\;\forall\;u\in E.
\end{equation}
We fix $\rho\in(0,1)$ such that $\rho<1/c_1$. Then  relation  \eq{el1} implies 
$$|u|_{q(x)}<1,\;\;\;\forall\;u\in E,\;{\rm with}\;\|u\|=\rho.$$
Furthermore, relation \eq{L5} yields
\begin{equation}\label{el2}
\int_\Omega|u|^{q(x)}\;dx\leq|u|_{q(x)}^{q^-},\;\;\;\forall\;u\in E,\;{\rm with}\;\|u\|=\rho.
\end{equation}
Relations \eq{el1} and \eq{el2} imply
\begin{equation}\label{el3}
\int_\Omega|u|^{q(x)}\;dx\leq c_1^{q^-}\|u\|^{q^-},\;\;\;\forall\;u\in E,\;{\rm with}\;\|u\|=\rho.
\end{equation}
Taking into account  relations \eq{L5} and \eq{el3} we deduce that 
for any $u\in E$ with $\|u\|=\rho$ the following inequalities hold true
\begin{eqnarray*}
J_\lambda(u)&\geq&\frac{1}{p^+}\int_{\Omega}|\nabla u|^{p(x)}\;dx-\frac{\lambda}{q^-}
\int_\Omega|u|^{q(x)}\;dx\\
&\geq&\frac{1}{p^+}\|u\|^{p^+}-\frac{\lambda}{q^-}c_1^{q^-}\|u\|^{q^-}\\
&=&\frac{1}{p^+}\rho^{p^+}-\frac{\lambda}{q^-}c_1^{q^-}\rho^{q^-}\\
&=&\rho^{q^-}\left(\frac{1}{p^+}\rho^{p^+-q^-}-\frac{\lambda}{q^-}c_1^{q^-}\right).
\end{eqnarray*}
By the above inequality we remark that if we define
\begin{equation}\label{el4}
\lambda^\star=\frac{\rho^{p^+-q^-}}{2p^+}\cdot\frac{q^-}{c_1^{q^-}}
\end{equation}
then for any $\lambda\in(0,\lambda^\star)$ and any $u\in E$ with $\|u\|=\rho$ there 
exists $a=\frac{\rho^{p^+}}{2p^+}>0$ such that
$$J_\lambda(u)\geq a>0.$$
The proof of Lemma \ref{l1} is complete.  \endproof

\begin{lemma}\label{l2}
There exists $\phi\in E$ such that $\phi\geq 0$,
$\varphi\neq 0$ and
$J_\lambda(t\phi)<0$,
for $t>0$ small enough. 
\end{lemma}

\proof
Assumption \eq{3} implies that $q^-<p^-$. Let $\epsilon_0>0$ be such that $q^-+\epsilon_0<p^-$. 
On the other hand, since $q\in C(\overline\Omega)$ it follows that there exists an open set 
$\Omega_0\subset\Omega$ such that $|q(x)-q^-|<\epsilon_0$ for all $x\in\Omega_0$. Thus, we 
conclude that $q(x)\leq q^-+\epsilon_0<p^-$ for all $x\in\Omega_0$.

Let $\phi\in C_0^\infty(\Omega)$ be such that ${\rm supp}(\phi)\supset\overline\Omega_0$, 
$\phi(x)=1$ for all $x\in\overline\Omega_0$ and $0\leq\phi\leq 1$ in $\Omega$. Then using the 
above information for any $t\in(0,1)$ we have
\begin{eqnarray*}
J_\lambda(t\phi)&=&\int_{\Omega}\frac{t^{p(x)}}{p(x)}|\nabla\phi|^{p(x)}\;dx-\lambda
\int_\Omega\frac{t^{q(x)}}{q(x)}|\phi|^{q(x)}\;dx\\
&\leq&\frac{t^{p^-}}{p^-}\int_{\Omega}|\nabla\phi|^{p(x)}\;dx-\frac{\lambda}{q^+}
\int_\Omega t^{q(x)}|\phi|^{q(x)}\;dx\\
&\leq&\frac{t^{p^-}}{p^-}\int_{\Omega}|\nabla\phi|^{p(x)}\;dx-\frac{\lambda}{q^+}
\int_{\Omega_0} t^{q(x)}|\phi|^{q(x)}\;dx\\
&\leq&\frac{t^{p^-}}{p^-}\int_{\Omega}|\nabla\phi|^{p(x)}\;dx-\frac{\lambda\cdot 
t^{q^-+\epsilon_0}}{q^+}\int_{\Omega_0}|\phi|^{q(x)}\;dx.
\end{eqnarray*}
Therefore
$$J_\lambda(t\phi)<0$$  
for $t<\delta^{1/(p^--q^--\epsilon_0)}$ with
$$0<\delta<\min\left\{1,\frac{\frac{\lambda\cdot p^-}{q^+}
\int_{\Omega_0}|\phi|^{q(x)}\;dx}
{\int_{\Omega}|\nabla\phi|^{p(x)}\;dx}\right\}.$$
Finally, we point out that $\int_{\Omega}|\nabla\phi|^{p(x)}\;dx>0$. Indeed, it is clear that 
$$\int_{\Omega_0}|\phi|^{q(x)}\;dx\leq\int_{\Omega}|\phi|^{q(x)}\;dx\leq
\int_{\Omega_0}|\phi|^{q^-}\;dx.$$
On the other hand, $W_0^{1,p(x)}(\Omega)$ is continuously embedded in $L^{q^-}(\Omega)$ and thus, 
there exists a positive constant $c_2$ such that
$$|\phi|_{q^-}\leq c_2\|\phi\|.$$
The last two inequalities imply that
$$\|\phi\|>0$$
and combining that fact with relations \eq{L4} or \eq{L5} we deduce that
$$\int_{\Omega}|\nabla\phi|^{p(x)}\;dx>0.$$
The proof of Lemma \ref{l2} is complete.   \endproof

\medskip
{\sc Proof of Theorem \ref{t1}.}
Let $\lambda^\star>0$ be defined as in \eq{el4} and $\lambda\in(0,\lambda^\star)$.
By Lemma \ref{l1} it follows that on the boundary of the  ball centered at the origin and of
radius $\rho$ in $E$, denoted by $B_\rho(0)$, we have
\begin{equation}\label{et8}
\inf\limits_{\partial B_\rho(0)}J_\lambda>0.
\end{equation}
On the other hand, by Lemma \ref{l2}, there exists $\phi\in E$ such that 
$J_\lambda(t\phi)<0$ for all $t>0$ small enough. Moreover, relations \eq{el3} and 
\eq{L5} imply that for any $u\in B_\rho(0)$ we have
$$J_\lambda(u)\geq\frac{1}{p^+}\|u\|^{p^+}-\frac{\lambda}{q^-}c_1^{q^-}\|u\|^{q^-}\,.$$
It follows that
$$-\infty<\underline{c}:=\inf\limits_{\overline{B_\rho(0)}}J_\lambda<0.$$
We let now $0<\epsilon<\inf_{\partial B_\rho(0)}J_\lambda-\inf_{B_\rho(0)}J_\lambda$. 
Applying Ekeland's variational principle to the functional
$J_\lambda:\overline{B_\rho(0)}\rightarrow\RR$,  we find 
$u_\epsilon\in\overline{B_\rho(0)}$ such that
\begin{eqnarray*}
J_\lambda(u_\epsilon)&<&\inf\limits_{\overline{B_\rho(0)}}J_\lambda+\epsilon\\
J_\lambda(u_\epsilon)&<&J_\lambda(u)+\epsilon\cdot\|u-u_\epsilon\|,\;\;\;u\neq u_\epsilon.
\end{eqnarray*} 
Since
$$J_\lambda(u_\epsilon)\leq\inf\limits_{\overline{B_\rho(0)}}J_\lambda+\epsilon\leq
\inf\limits_{B_\rho(0)} J_\lambda+\epsilon<\inf\limits_{\partial B_\rho(0)}J_\lambda\,,$$
we deduce that $u_\epsilon\in B_\rho(0)$. Now, we define $I_\lambda:
\overline{B_\rho(0)}\rightarrow\RR$ by $I_\lambda(u)=J_\lambda(u)+\epsilon\cdot
\|u-u_\epsilon\|$. It is clear that $u_\epsilon$ is a minimum point 
of $I_\lambda$ and thus
$$\di\frac{I_\lambda(u_\epsilon+t\cdot v)-{I_\lambda}(u_\epsilon)}
{t}\geq 0$$
for small $t>0$ and any $v\in B_1(0)$. The
above relation yields
$$\di\frac{J_\lambda(u_\epsilon+t\cdot v)-J_\lambda(u_\epsilon)}{t}+
\epsilon\cdot\|v\|\geq 0.$$
Letting $t\rightarrow 0$ it follows that $\langle J_\lambda^{'}
(u_\epsilon),v\rangle+\epsilon\cdot\|v\|>0$ and we infer that 
$\|J_\lambda^{'}(u_\epsilon)\|\leq\epsilon$.

We deduce that there exists a sequence
$\{w_n\}\subset B_\rho(0)$ such that 
\begin{equation}\label{et9}
J_\lambda(w_n)\rightarrow{\underline c}\;\;\;{\rm and}\;\;\;
J_\lambda^{'}(w_n)\rightarrow 0.
\end{equation}
It is clear that $\{w_n\}$ is bounded in $E$. Thus, there exists $w\in E$ such that, up  to a
subsequence, $\{w_n\}$ converges weakly to $w$ in $E$. 
Since $q(x)<p^\star(x)$ for all $x\in\overline\Omega$ we deduce that $E$ is compactly embedded in 
$L^{q(x)}(\Omega)$, hence $\{w_n\}$ converges strongly to $w$ in $L^{q(x)}(\Omega)$. 
So, by relations \eq{L6} and \eq{Hol},
$$\lim\limits_{n\rightarrow\infty}\int_\Omega|w_n|^{q(x)-2}w_n(w_n-w)\;dx=0.$$
On the other hand, relation \eq{et9} yields
$$\lim\limits_{n\rightarrow\infty}\langle J_\lambda^{'}(w_n),w_n-w\rangle=0.$$
Using the above information we find
\begin{equation}\label{et10}
\lim\limits_{n\rightarrow\infty}\int_\Omega|\nabla w_n|^{p(x)-2}\nabla w_n
\nabla(w_n-w)\;dx=0.
\end{equation}
Relation \eq{et10} and the fact that $\{w_n\}$ converges weakly to $w$ in $E$ 
enable us to apply Theorem 3.1 in Fan and Zhang \cite{FZh} in order to obtain that $\{w_n\}$ 
converges strongly to $w$ in $E$. So, by \eq{et9},
\begin{equation}\label{et11}
J_\lambda(w)=\underline c<0\;\;\;{\rm and}\;\;\;J_\lambda^{'}(w)=0.
\end{equation}
We conclude that $w$ is a nontrivial weak solution for problem \eq{2} and thus any 
$\lambda\in(0,\lambda^\star)$ is an eigenvalue of problem \eq{2}.

The proof of Theorem \ref{t1} is complete.\endproof

\medskip
Let us now assume that the hypotheses of Theorem \ref{t1} are fulfilled and,
furthermore,
$$\max_{\overline\Omega}p(x)<\max_{\overline\Omega}q(x).$$
Then, using similar arguments as in the proof of Lemma \ref{l2}, we find some 
$\psi\in E$ such that
$$\lim\limits_{t\rightarrow\infty}J_\lambda(t\psi)=-\infty.$$
That fact combined with Lemma \ref{l1} and the mountain pass theorem (see \cite{AR}) implies that there 
exists a sequence $\{u_n\}$ in $E$ such that
\begin{equation}\label{mp}
J_\lambda(u_n)\rightarrow \overline{c}>0\;\;\;{\rm and}\;\;\;
J_\lambda^{'}(u_n)\rightarrow 0\;{\rm in}\; E^\star.
\end{equation}
However, relation \eq{mp} is not useful because we can not show that the sequence $\{u_n\}$ 
is bounded in $E$ since the functional $J_\lambda$ does not satisfy a relation of the
Ambrosetti-Rabinowitz type. This enable us to affirm that we can not obtain a critical 
point for $J_\lambda$ by using this method.

On the other hand, we point out that we will fail in trying to show that the functional 
$J_\lambda$ is coercive since by relation \eq{3} we have $q^+>p^-$. Thus, we can not apply
(as in the homogeneous case) a result 
as Theorem 1.2 in Struwe \cite{S} in order to obtain a critical point of the functional 
$J_\lambda$.


\begin{thebibliography}{99}
{\footnotesize

\bibitem{AR} A. Ambrosetti and P. H. Rabinowitz, Dual variational
methods in critical point theory,
{\it J. Funct. Anal.} {\bf 14} (1973), 349-381.

\bibitem{anane} A. Anane,  Simplicit\'e et isolation de la premi\`ere valeur propre du $p$-laplacien avec poids, 
{\it  C. R. Acad. Sci. Paris S\'er. I} {\bf 305}  (1987),  725-728.

\bibitem{boc} M. Bocher, The smallest characteristic numbers in a certain exceptional case,
{\it Bull. Amer. Math. Soc.} {\bf 21} (1914), 6-9.

\bibitem{D} L. Diening, {\it Theoretical and Numerical Results for
Electrorheological Fluids}, Ph.D. thesis, University of Frieburg,
Germany, 2002.

\bibitem{edm} D. E. Edmunds, J. Lang, and A. Nekvinda, On $L^{p(x)}$
norms, {\it Proc. Roy. Soc. London Ser.~A} {\bf 455} (1999), 219-225.

\bibitem{edm2} D. E. Edmunds and J. R\'akosn\'{\i}k, Density of smooth
functions in $W^{k,p(x)}(\Omega)$, {\it Proc. Roy. Soc. London Ser.~A}
{\bf 437} (1992), 229-236.

\bibitem{edm3} D. E. Edmunds and J. R\'akosn\'{\i}k, Sobolev embedding
with variable exponent, {\it Studia Math.} {\bf 143} (2000), 267-293.

\bibitem{E} I. Ekeland, On the variational principle, {\it J. Math.
Anal. Appl.} {\bf 47} (1974), 324-353.

\bibitem{FZh} X. L. Fan and Q. H. Zhang, Existence of solutions for
$p(x)$-Laplacian Dirichlet problem, {\it Nonlinear Anal.} {\bf 52}
(2003), 1843-1852.

\bibitem{FZZ} X. Fan, Q. Zhang and D. Zhao, Eigenvalues of $p(x)$-Laplacian 
Dirichlet problem, {\it J. Math. Anal. Appl.} {\bf 302} (2005), 306-317.

\bibitem{hal} T. C. Halsey, Electrorheological fluids, {\it Science} 
{\bf 258} (1992), 761-766.

\bibitem{hes} P. Hess and T. Kato, On some linear and nonlinear eigenvalue problems with an indefinite weight function, 
{\it Comm. Partial Differential Equations} {\bf 5}  (1980),  999-1030. 

\bibitem{KR} O. Kov\'a\v cik and J. R\'akosn\'{\i}k, On spaces 
$L^{p(x)}$ and
$W^{1,p(x)}$, {\it Czechoslovak Math. J.} {\bf 41} (1991), 592-618.

\bibitem{RoyalSoc} M. Mih\u ailescu and V. R\u adulescu, 
A multiplicity result for a nonlinear degenerate problem
arising in the theory of electrorheological fluids, {\it Proc. Roy. Soc. London Ser.~A}, in press.

\bibitem{minple} S. Minakshisundaram and A. Pleijel, Some properties of the eigenfunctions of the Laplace-operator on Riemannian manifolds, {\it Canadian J. Math.} {\bf  1}  (1949), 242-256. 

\bibitem{M} J. Musielak, {\it Orlicz Spaces and Modular  Spaces},
Lecture Notes in Mathematics, Vol. 1034, Springer, Berlin, 1983.

\bibitem{ple} A. Pleijel, On the eigenvalues and eigenfunctions of elastic plates, 
{\it Comm. Pure Appl. Math.}  {\bf 3}  (1950), 1-10.

\bibitem{R} M. Ruzicka, {\it Electrorheological Fluids: Modeling
and Mathematical Theory}, Springer-Verlag, Berlin, 2002.

\bibitem{samko} S. Samko and B. Vakulov, Weighted Sobolev theorem with variable exponent
for spatial and spherical potential operators, {\it J. Math. Anal. Appl.} {\bf 310} (2005), 229-246.

\bibitem{S} M. Struwe, {\it Variational Methods: Applications to
Nonlinear Partial Differential Equations and Hamiltonian Systems},
Springer, Heidelberg, 1996.

\bibitem{Z1} V. Zhikov, Averaging of functionals in the calculus of
variations and elasticity, {\it Math. USSR Izv.} {\bf 29} (1987), 
33-66.

}
\end{thebibliography}
\end{document}